\documentclass[11pt]{article}

\usepackage[latin1]{inputenc} 
\usepackage{amsthm,amsmath,amssymb} 
\usepackage{graphicx}
\usepackage{color}
\usepackage{verbatim}

\usepackage{enumitem} 
\setlist[enumerate]{topsep=0pt,itemsep=-1ex,partopsep=1ex,parsep=1ex}

\theoremstyle{plain}

\theoremstyle{definition}

\newcommand{\mb}[1]{\mathbb{#1}}
\newcommand{\nib}[1]{\noindent {\bf #1}}

\newcommand{\sub}{\subseteq}

\newcommand{\sups}{\supseteq}

\newcommand{\sm}{\setminus}

\newcommand{\dD}{\delta}

\newcommand{\lL}{\lambda}
\newcommand{\tT}{\theta}

\newcommand{\OO}{\Omega}

\title{Hypergraph matchings and designs}

\author{Peter Keevash\thanks{Mathematical Institute,
University of Oxford, Oxford, UK. Email: keevash@maths.ox.ac.uk.
\newline \hspace*{1.8em}Research supported
in part by ERC Consolidator Grant 647678.}}

\begin{document}

\maketitle

\begin{abstract}
We survey some aspects of the perfect matching problem in hypergraphs,
with particular emphasis on structural characterisation of the existence 
problem in dense hypergraphs and the existence of designs.
\end{abstract}

\section{Introduction}

Matching theory is a rich and rapidly developing subject
that touches on many areas of Mathematics and its applications.
Its roots are in the work of Steinitz \cite{St}, Egerv\'ary \cite{eg},
Hall \cite{hall} and K\"onig \cite{konig} 
on conditions for matchings in bipartite graphs.
After the rise of academic interest in efficient algorithms
during the mid 20th century, three cornerstones of matching theory were
Kuhn's \cite{kuhn} `Hungarian' algorithm for the Assignment Problem,
Edmonds' \cite{edmonds} algorithm for finding a maximum matching
in general (not necessarily bipartite) graphs,
and the Gale-Shapley \cite{GS} algorithm for Stable Marriages.
For an introduction to matching theory in graphs 
we refer to Lov\'asz and Plummer \cite{LP},
and for algorithmic aspects to parts II and III
of Schrijver \cite{schrijver}.

There is also a very large literature on matchings 
in hypergraphs. This article will be mostly concerned 
with one general direction in this subject, 
namely to determine conditions under which the
necessary `geometric' conditions of `space' and `divisibility'
are sufficient for the existence of a perfect matching.
We will explain these terms and discuss some aspects of this
question in the next two sections, but first, for the remainder 
of this introduction, we will provide some brief pointers 
to the literature in some other directions.

We do not expect a simple general characterisation of the
perfect matching problem in hypergraphs, 
as by contrast with the graph case,
it is known to be NP-complete even for $3$-graphs
(i.e.\ when all edges have size $3$),
indeed, this was one of Karp's original 
21 NP-complete problems~\cite{Karp}.
Thus for algorithmic questions related to hypergraph matching,
we do not expect optimal solutions, and may instead consider 
Approximation Algorithms (see e.g.\ \cite{WS, AFS, LRS}).

Another natural direction is to seek nice sufficient conditions
for perfect matchings. There is a large literature in Extremal
Combinatorics on results under minimum degree assumptions,
known as `Dirac-type' theorems, after the classical result
of Dirac \cite{dirac} that any graph on $n \ge 3$ vertices
with minimum degree at least $n/2$ has a Hamiltonian cycle.
It is easy to see that $n/2$ is also the minimum degree threshold
for a graph on $n$ vertices (with $n$ even) to have a perfect matching,
and this exemplifies the considerable similarities between the
perfect matching and Hamiltonian problems (but there are also
substantial differences). A landmark in the efforts to obtain
hypergraph generalisations of Dirac's theorem was the result
of R\"odl, Ruci\'nski and Szemer\'edi \cite{RRS} that determined
the codegree threshold for perfect matchings in uniform hypergraphs;
this paper was significant for its proof method as well as the result,
as it introduced the Absorbing Method (see section \ref{sec:ab}), 
which is now a very important tool 
for proving the existence of spanning structures.
There is such a large body of work in this direction
that it needs several surveys to describe, and indeed
these surveys already exist \cite{RR, KOlarge, KOicm, yi, Y}.
The most fundamental open problem in this area 
is the Erd\H{o}s Matching Conjecture \cite{ematch},
namely that the maximum number of edges in an $r$-graph\footnote{
An $r$-graph is a hypergraph in which every edge contains $r$ vertices.}
on $n$ vertices with no matching of size $t$ is either
achieved by a clique of size $tr-1$ or the set of all
edges hitting some fixed set of size $t-1$
(see \cite[Section 3]{FT} for discussion and
a summary of progress).

The duality between matching and covers in hypergraphs
is of fundamental important in Combinatorics (see \cite{fur})
and Combinatorial Optimisation (see \cite{C}).
A defining problem for this direction of research
within Combinatorics is `Ryser's Conjecture'
(published independently by Henderson \cite{hen}
and Lov\'asz \cite{L}) that in any $r$-partite $r$-graph
the ratio of the covering and matching numbers is at most $r-1$.
For $r=2$ this is K\"onig's Theorem. The only other known
case is $r=3$, due to Aharoni \cite{A}, using a hypergraph
analogue of Hall's theorem due to Aharoni and Haxell \cite{AH},
which has a topological proof. There are now many applications
of topology to hypergraph matching, and more generally
`independent transversals' (see the survey \cite{haxsurvey}).
In the other direction, the hypergraph matching complex
is now a fundamental object of Combinatorial Topology,
with applications to Quillen complexes in Group Theory,
Vassiliev knot invariants and Computational Geometry
(see the survey \cite{wachs}).

From the probabilistic viewpoint, 
there are (at least) two natural questions: 

(i) does a random hypergraph
have a perfect matching with high probability (whp)?

(ii) what does a random matching 
from a given (hyper)graph look like?

The first question for the usual (binomial) random hypergraph
was a longstanding open problem, perhaps first stated by
Erd\H{o}s \cite{e81} (who attributed it to Shamir),
finally solved by Johansson, Kahn and Vu \cite{JKV};
roughly speaking, the threshold is `where it should be',
namely around the edge probability at which with high
probability every vertex is in at least one edge.
Another such result due to
Cooper, Frieze, Molloy and Reed \cite{CFMR}
is that random regular hypergraphs
(of fixed degree and edge size)
whp have perfect matchings.

The properties of random matchings in lattices have been 
extensively studied under the umbrella of the `dimer model' 
(see \cite{kenyon}) in Statistical Physics. 
However, rather little is known regarding the typical structure
of random matchings in general graphs, let alone hypergraphs.
Substantial steps in this direction have been taken by results
of Kahn \cite{Ka2} characterising when the size of 
a random matching has an approximate normal distribution,
and Kahn and Kayll \cite{KaKa} establishing long-range decay
of correlations of edges in random matchings in graphs;
the final section of \cite{Ka2} contains many open problems,
including conjectural extensions to simple hypergraphs.

Prequisite to the understanding of random matchings
are the closely related questions of Sampling and Approximate Counting
(as established in the Markov Chain Monte Carlo 
framework of Jerrum and Sinclair, see \cite{J}).
An approximate counting result for hypergraph matchings 
with respect to balanced weight functions was obtained 
by Barvinok and Samorodnitsky \cite{BS}.
Extremal problems also arise naturally in this context, 
for the number of matchings, and more generally for
other models in Statistical Physics, such as the 
hardcore model for independent sets.
Much of the recent progress here appears in the survey \cite{yufei},
except for the very recent solution of (almost all cases of)
the Upper Matching Conjecture of
Friedland, Krop and Markstr\"om \cite{FKM}
by Davies, Jenssen, Perkins and Roberts \cite{DJPR}.

\section{Space and divisibility}

In this section we discuss a result
(joint work with Mycroft \cite{KM})
that characterises the obstructions to perfect matching
in dense hypergraphs (under certain conditions
to be discussed below). The obstructions are geometric
in nature and are of two types:
Space Barriers (metric obstructions)
and Divisibility Barriers (arithmetic obstructions).

The simplest illustration of these two phenomena
is seen by considering extremal examples for the
simple observation mentioned earlier that a graph
on $n$ vertices ($n$ even) with minimum degree 
at least $n/2$ has a perfect matching.
One example of a graph with minimum degree $n/2-1$
and no perfect matching is obtained by fixing a set
$S$ of $n/2-1$ vertices and taking all edges that
intersect $S$. Then in any matching $M$, each edge
of $M$ uses at least one vertex of $S$, 
so $|M| \le |S| < n/2$; there is no `space'
for a perfect matching. For another example,
suppose $n=2$ mod $4$ and consider the graph
that is the disjoint union of two cliques each
of size $n/2$ (which is odd). As edges have size $2$,
which is even, there is an arithmetic (parity)
obstruction to a perfect matching.

There is an analogous parity obstruction to matching
in general $r$-graphs, namely an $r$-graph $G$
with vertices partitioned as $(A,B)$, so that 
$|A|$ is odd and $|e \cap A|$ is even
for each edge $e$ of $G$; this is one of
the extremal examples for the codegree threshold
of perfect matchings (see \cite{RRS}). 

In general, space barriers are constructions 
for each $1 \le i \le r$, obtained by fixing a set $S$ 
of size less than $in/r$ and taking the $r$-graph of 
all edges $e$ with $|e \cap S| \ge i$. Then for any matching 
$M$ we have $|M| \le |S|/i < n/r$, so $M$ is not perfect.

General divisibility barriers are obtained
by fixing a lattice (additive subgroup) $L$ in $\mb{Z}^d$
for some $d$, fixing a vertex set partitioned
as $(V_1,\dots,V_d)$, with $(|V_1|,\dots,|V_d|) \notin L$,
and taking the $r$-graph of all edges $e$ such that 
$(|e \cap V_1|,\dots,|e \cap V_d|) \in L$.
For example, the parity obstruction corresponds 
to the lattice $\{(2x,y): x,y \in \mb{Z}\}$.

To state the result of \cite{KM} that is most conveniently applicable
we introduce the setting of simplicial complexes and
degree sequences. We consider a simplicial complex $J$
on $[n]=\{1,\dots,n\}$, write $J_i = \{e \in J: |e|=i\}$ and look
for a perfect matching in the $r$-graph $J_r$.
We define the degree sequence $(\dD_0(J),\dots,\dD_{r-1}(J))$
so that each $\dD_i(J)$ is the least $m$ such that
each $e \in J_i$ is contained in at least $m$
edges of $J_{i+1}$. We define the critical degree sequence 
$\dD^c = (\dD^c_0,\dots,\dD^c_{r-1})$ by $\dD^c_i=(1-i/r)n$.
The space barrier constructions show that for each $i$
there is a complex with $\dD_i(J)$ slightly less than
$\dD^c_i$ but no perfect matching.
An informal statement of \cite[Theorem 2.9]{KM}
is that if $J$ is an $r$-complex on $[n]$ 
(where $r \mid n$) with all $\dD_i(J) \ge \dD_i^c - o(n)$
such that $J_r$ has no perfect matching
then $J$ is close (in edit distance)
to a space barrier or a divisibility barrier.

One application of this result (also given in \cite{KM})
is to determine the exact codegree threshold for packing
tetrahedra in $3$-graphs; it was surprising that it was
possible to obtain such a result given that the simpler-sounding
problems of determining the thresholds (edge or codegree)
for the existence of just one tetrahedron are open,
even asymptotically (the edge threshrold is a famous
conjecture of Tur\'an; for more on Tur\'an problems
for hypergraphs see the survey \cite{Kturan}).
Other applications are a multipartite version of the
Hajnal-Szemeredi theorem (see \cite{KM2})
and determining the `hardness threshold'
for perfect matchings in dense hypergraphs
(see \cite{KKM, han}).
 
We will describe the hardness threshold in more detail,
as it illustrates some important features of space
and divisibility, and the distinction between
perfect matchings and almost perfect matchings.
For graphs there is no significant difference
in the thresholds for these problems,
whereas for general $r$-graphs there is 
a remarkable contrast: the codegree threshold
for perfect matchings \cite{RRS} is about $n/2$,
whereas Han \cite{han2}, proving a conjecture from \cite{RRS},
showed that a minimum codegree of only $n/r$ guarantees
a matching of size $n/r-1$, i.e.\ one less than perfect.
The explanation for this contrast is that the divisibility
barrier is no obstacle to almost perfect matching,
whereas the space barrier is more robust,
and can be continuously `tuned' to 
exclude a matching of specified size.

To illustrate this, we consider a $3$-graph $G_0$ on $[n]$ where 
the edges are all triples that intersect some fixed set $S$
of size $(1/3 - c)n$, for some small $c>0$.
Then the minimum codegree and maximum matching size
in $G_0$ are both equal to $|S|$. Furthermore,
if we consider $G = G_0 \cup G_1$ where all edges
of $G_1$ lie within some $S'$ disjoint from $S$
with $|S'|=3cn$ then $G$ has a perfect matching
if and only if $G_1$ has a perfect matching,
which is NP-complete to decide for arbitrary $G_1$.
Thus the robustness of the space barrier provides
a reduction showing that the codegree threshold
for the existence of an algorithm for the perfect
matching is at least the threshold
for an approximate perfect matching.

Now consider the decision problem for perfect matchings 
in $3$-graphs on $[n]$ (where $3 \mid n$) with 
minimum codegree at least $\dD n$. For $\dD<1/3$
the problem is NP-complete, and for $\dD>1/2$
it is trivial (there is a perfect matching by \cite{RRS}).
For intermediate $\dD$ there 
is a polynomial-time algorithm, and this is in essence
a structural stability result: the main ingredient of
the algorithm is a result of \cite{KKM} that any such
$3$-graph with no perfect matching is contained 
in a divisibility barrier. (For general $r$ the 
structural characterisation is more complicated.)

\section{Fractional matchings}

The key idea of the Absorbing Method \cite{RRS}
mentioned earlier is that the task of finding
perfect matchings can often broken into two subproblems:
(i) finding almost perfect matchings,
(ii) absorbing uncovered vertices into 
an almost perfect matching until it becomes perfect.
We have already seen that the almost perfect matching
problem appears naturally as a relaxation of the perfect
matching problem in which we eliminate divisibility
obstacles but retain space obstacles. This turns out
to fit into a more general framework of fractional matchings,
in which the relaxed problem is a question of convex geometry,
and space barriers correspond to separating hyperplanes.

The fractional (linear programming) relaxation of the perfect matching
problem in a hypergraph is to assign non-negative weights
to the edges so that for any vertex $v$, there is a total 
weight of $1$ on all edges incident to $v$. A perfect matching 
corresponds to a $\{0,1\}$-valued solution, so the existence
of a fractional perfect matching is necessary for the existence
of a perfect matching. We can adopt a similar point of view
regarding divisibility conditions. Indeed, we can similarly define
the integer relaxation of the perfect matching problem in which we now
require the weights to be integers (not necessarily non-negative);
then the existence of an integral perfect matching is 
necessary for the existence of a perfect matching.

The fractional matching problem appears naturally in Combinatorial Optimisation
(see \cite{C, schrijver}) because it brings in
polyhedral methods and duality to bear on the matching problem.
It has also been studied as a problem in its own right from the
perspective of random thresholds (e.g.\ \cite{DK, kriv}),
and it appears naturally in combinatorial existence problems,
as in dense hypergraphs almost perfect matchings and fractional matchings 
tend to appear at the same threshold. Indeed, for many
open problems, such as the Erd\H{o}s Matching Conjecture \cite{ematch}
or the Nash-Williams Triangle Decomposition Conjecture \cite{NW},
any progress on the fractional problem translates
directly into progress on the original problem (see \cite{BKLO}).

This therefore makes the threshold problem for fractional matchings
and decompositions a natural problem in its own right.
For example, an asymtotic solution of the Nash-Williams Conjecture
would follow from the following conjecture: any graph on $n$ vertices
with minimum degree at least $3n/4$ has a fractional triangle decomposition,
i.e.\ an assignment of non-negative weights to its triangles
so that for any edge $e$ there is total weight $1$
on the triangles containing $e$. An extremal example $G$
for this question can be obtained by taking a balanced complete bipartite
graph $H$ and adding a $(n/4-1)$-regular graph inside each part;
indeed, this is a space barrier to a fractional triangle decomposition,
as any triangle uses at least one edge not in $H$, but $|H| > 2|G \sm H|$.
The best known upper bound is $0.913n$ by Dross \cite{dross}.
More generally, Barber, K\"uhn, Lo, Montgomery and Osthus \cite{BKLMO}
give the current best known bounds on the thresholds for fractional
clique decompositions (in graphs and hypergraphs),
but these seem to be far from optimal.

There are (at least) two ways to think about the relationship
between almost perfect matchings and fractional matchings.
The first goes back to the `nibble' (semi-random) method of R\"odl \cite{R},
introduced to solve the Erd\H{o}s-Hanani \cite{EH} conjecture on
approximate Steiner systems (see the next section),
which has since had a great impact on Combinatorics 
(e.g.\ \cite{AKS,BB,Boh,BFL,BK,BK2,FGM,FR,Gr,Ka,KaLP,Kim,KR,Kuz,PS,S,Vu}).
A special case of a theorem of Kahn \cite{KaLP}
is that if there is a fractional perfect matching
on the edges of an $r$-graph $G$ on $[n]$ 
such that for any pair of vertices $x,y$
the total weight on edges containing $\{x,y\}$
is $o(1)$ then $G$ has a matching covering
all but $o(n)$ vertices. In this viewpoint,
it is natural to interpret the weights of a fractional
matching as probabilities, and an almost perfect matching
as a random rounding; in fact, this random rounding 
is obtained iteratively, so there are some parallels with 
the development of iterative rounding algorithms (see \cite{LRS}).

Another way to establish the connection between
almost perfect matchings and fractional matchings
is via the theory of Regularity, developed by
Szemer\'edi \cite{Sz} for graphs and extended to hypergraphs
independently by Gowers \cite{Gow} and
Nagle, R\"odl, Schacht and Skokan \cite{NRS, RSc1, RSc2, RSk}.
(The connection was first established by
Haxell and R\"odl \cite{HR} for graphs and
R\"odl, Schacht, Siggers and Tokushige \cite{RSST} for hypergraphs.)
To apply Regularity to obtain spanning structures
(such as perfect matchings) requires an accompanying result
known as a blowup lemma, after the original such result for graphs 
obtained by Koml\'os, S\'ark\"ozy and Szemer\'edi \cite{KSS};
we proved the hypergraph version in \cite{Kblowup}.
More recent developments (for graphs) along these lines 
include the Sparse Blowup Lemmas of 
Allen, B\"ottcher, H\`an,  Kohayakawa and Person \cite{ABHKP}
and a blowup-up lemma suitable for decompositions
(as in the next section) obtained by
Kim, K\"uhn, Osthus and Tyomkyn \cite{KKOT}
(it would be interesting and valuable 
to obtain hypergraph versions of these results).
The technical difficulties of the Hypergraph Regularity method 
are a considerable barrier to its widespread application,
and preclude us giving here a precise statement of 
\cite[Theorem 9.1]{KM}, which informally speaking
characterises the perfect matching problem in dense
hypergraphs with certain extendability conditions
in terms of space and divisibility.

\section{Designs and decompositions}

A \emph{Steiner system} with parameters $(n,q,r)$ 
is a $q$-graph $G$ on $[n]$ such that any $r$-set 
of vertices is contained in exactly one edge.
For example, a Steiner Triple System on $n$ points
has parameters $(n,3,2)$.
The question of whether there is a Steiner system with given 
parameters is one of the oldest problems in combinatorics, 
dating back to work of Pl\"ucker (1835), Kirkman (1846) 
and Steiner (1853); see \cite{RobinW} for a historical account. 

Note that a Steiner system with parameters $(n,q,r)$
is equivalent to a $K^r_q$-decomposition of $K^r_n$
(the complete $r$-graph on $[n]$).
It is also equivalent to a perfect matching
in the auxiliary $\tbinom{q}{r}$-graph on $\tbinom{[n]}{r}$
(the $r$-subsets of $[n]:=\{1,\dots,n\}$) with edge set
$\{ \tbinom{Q}{r}: Q \in \tbinom{[n]}{q} \}$.

More generally, we say that a set $S$ of $q$-subsets of an 
$n$-set $X$ is a \emph{design} with parameters $(n,q,r,\lL)$ if 
every $r$-subset of $X$ belongs to exactly $\lL$ elements of $S$.
(This is often called an `$r$-design' in the literature.)
There are some obvious necessary `divisibility conditions' 
for the existence of such $S$, namely that $\tbinom{q-i}{r-i}$ 
divides $\lL \tbinom{n-i}{r-i}$ for every $0 \le i \le r-1$ 
(fix any $i$-subset $I$ of $X$ and consider the sets in $S$ 
that contain $I$). It is not known who first advanced the 
`Existence Conjecture' that the divisibility conditions are 
also sufficient, apart from a finite number of exceptional $n$ 
given fixed $q$, $r$ and $\lL$. 

The case $r=2$ has received particular attention due to 
its connections to statistics, under the name of 
`balanced incomplete block designs'. 
We refer the reader to \cite{CD} for a summary
of the large literature and applications of this field. 
The Existence Conjecture for $r=2$ was a long-standing
open problem, eventually resolved by Wilson \cite{W1,W2,W3} 
in a series of papers that revolutionised Design Theory.
The next significant progress on the general conjecture
was in the solution of the two relaxations (fractional and integer)
discussed in the previous section 
(both of which are interesting in their own right 
and useful for the original problem).
We have already mentioned R\"odl's solution of
the Erd\H{o}s-Hanani Conjecture on approximate Steiner systems.
The integer relaxation was solved independently by
Graver and Jurkat \cite{GJ} and Wilson \cite{W4},
who showed that the divisibility conditions 
suffice for the existence of integral designs 
(this is used in \cite{W4} to show the existence 
for large $\lL$ of integral designs with non-negative coefficients).
Wilson \cite{W5} also characterised the existence 
of integral $H$-decompositions for any $r$-graph $H$. 

The existence of designs with $r \ge 7$ and any 
`non-trivial' $\lL$ was open before the breakthrough result 
of Teirlinck \cite{T} confirming this.
An improved bound on $\lL$ and a probabilistic method 
(a local limit theorem for certain random walks in high dimensions)
for constructing many other rigid combinatorial structures 
was recently given by Kuperberg, Lovett and Peled \cite{KLP}. 
Ferber, Hod, Krivelevich and Sudakov \cite{FHKS}
gave a construction of `almost Steiner systems',
in which every $r$-subset is covered by either one or two $q$-subsets.

In \cite{Kexist} we proved the Existence Conjecture in general, 
via a new method of Randomised Algebraic Constructions.
Moreover, in \cite{Kcount} we obtained the following
estimate for the number $D(n,q,r,\lL)$ of designs 
with parameters $(n,q,r,\lL)$ satisfying the necessary
divisibility conditions:
writing $Q=\tbinom{q}{r}$ and $N=\tbinom{n-r}{q-r}$, 
we have \[D(n,q,r,\lL) = \lL!^{-\tbinom{n}{r}} 
( (\lL/e)^{Q-1} N + o(N))^{\lL Q^{-1} \tbinom{n}{r}}.\]
Our counting result is complementary to that in \cite{KLP},
as it applies (e.g.) to Steiner Systems, whereas theirs
is only applicable to large multiplicities (but also allows
the parameters $q$ and $r$ to grow with $n$, and gives
an asymptotic formula when applicable).

The upper bound on the number of designs follows from
the entropy method pioneered by Radhakrishnan \cite{rad};
more generally, Luria \cite{luria} has recently established 
a similar upper bound on the number of perfect matchings
in any regular uniform hypergraph with small codegrees.
The lower bound essentially matches the number of choices
in the Random Greedy Hypergraph Matching process (see \cite{BB})
in the auxiliary $Q$-graph defined above,
so the key to the proof is showing that this process 
can be stopped so that whp it is possible to complete
the partial matching thus obtained to a perfect matching.
In other words, instead of a design, which can be viewed
as a $K^r_q$-decomposition of the $r$-multigraph $\lL K^r_n$,
we require a $K^r_q$-decomposition of some sparse submultigraph,
that satisfies the necessary divisibility conditions,
and has certain pseudorandomness properties
(guaranteed whp by the random process).

The main result of \cite{Kexist} achieved this, and indeed 
(in the second version of the paper) we obtained a more general 
result in the same spirit as \cite{KM}, namely that we can find
a clique decomposition of any $r$-multigraph
with a certain `extendability' property
that satisfies the divisibility conditions
and has a `suitably robust' fractional clique decomposition.

Glock, K\"uhn, Lo and Osthus \cite{GKLO,GKLO2}
have recently given a new proof of the existence of designs,
as well as some generalisations, including 
the existence of $H$-decompositions for any hypergraph $H$
(a question from \cite{Kexist}), relaxing the quasirandomness
condition from \cite{Kexist} (version 1) to an extendability
condition in the same spirit as \cite{Kexist} (version 2), 
and a more effective bound than that in \cite{Kexist}
on the minimum codegree decomposition threshold;
the main difference in our approaches lies in the treatment
of absorption (see the next section).

\section{Absorption} \label{sec:ab}

Over the next three sections we will sketch some approaches 
to what is often the most difficult part of a hypergraph
matching or decomposition problem, namely converting
an approximate solution into an exact solution.
We start by illustrating the Absorbing Method in its 
original form, namely the determination in \cite{RRS} of 
the codegree threshold for perfect matchings in $r$-graphs;
for simplicity we consider $r=3$.
 
We start by solving the almost perfect matching problem.
Let $G$ be a $3$-graph on $[n]$ with $3 \mid n$
and minimum codegree $\dD(G)=n/3$,
i.e.\ every pair of vertices is in at least $n/3$ edges.
We show that $G$ has a matching of size $n/3-1$ 
(i.e.\ one less than perfect). 
To see this, consider a maximum size matching $M$,
let $V_0 = V(G) \sm V(M)$, and suppose $|V_0|>3$.
Then $|V_0| \ge 6$, so we can fix disjoint
pairs $a_1b_1$, $a_2b_2$, $a_3b_3$ in $V_0$.
For each $i$ there are at least $n/3$ choices of $c$
such that $a_ib_ic \in E(G)$, and by maximality of $M$
any such $c$ lies in $V(M)$. We define the weight $w_e$
of each $e \in M$ as the number of edges of $G$ of the
form $a_ib_ic$ with $c \in e$.
Then $\sum_{e \in M} w_e \ge n$, and $|M|<n/3$,
so there is $e \in M$ with $w_e \ge 4$.
Then there must be distinct $c,c'$ in $e$
and distinct $i,i'$ in $[3]$ such that
$a_ib_ic$ and $a_{i'}b_{i'}c'$ are edges.
However, deleting $e$ and adding these edges
contradicts maximality of $M$.

Now suppose $\dD(G)=n/2 + cn$, 
where $c>0$ and $n>n_0(c)$ is large.
Our plan for finding a perfect matching
is to first put aside an `absorber' $A$,
which will be a matching in $G$ with the 
property that for any triple $T$ in $V(G)$
there is some edge $e \in A$ such that $T \cup e$
can be expressed as the disjoint union of two edges in $G$
(then we say that $e$ absorbs $T$).
Suppose that we can find such $A$,
say with $|A| < n/20$. Deleting the vertices of $A$
leaves a $3$-graph $G'$ on $n'=n-|A|$ vertices
with $\dD(G') \ge \dD(G)-3|A| > n'/3$.
As shown above, $G'$ has a matching $M'$ 
with $|M'| = n'/3-1$. Let $T = V(G') \sm V(M')$.
By choice of $A$ there is $e \in A$ such that
$T \cup e = e_1 \cup e_2$ for some 
disjoint edges $e_1,e_2$ in $G$.
Then $M' \cup (A \sm \{e\}) \cup \{e_1,e_2\}$
is a perfect matching in $G$.

It remains to find $A$. The key idea
is that for any triple $T$ there are many edges in $G$
that absorb $T$, and so if $A$ is random then
whp many of them will be present. We can bound 
the number of absorbers for any triple $T=xyz$
by choosing vertices sequentially. Say we want
to choose an edge $e=x'y'z'$ so that
$x'yz$ and $xy'z'$ are also edges.
There are at least $n/2 + cn$ choices for $x'$
so that $x'yz$ is an edge. Then for each
of the $n-4$ choices of $y' \in V(G) \sm \{x,y,z,x'\}$
there are at least $2cn - 1$ choices for
$z' \ne z$ so that $x'y'z'$ and $xy'z'$ are edges.
Multiplying the choices we see 
that $T$ has at least $cn^3$ absorbers.

Now suppose that we construct $A$ by choosing
each edge of $G$ independently with probability
$c/(4n^2)$ and deleting any pair that intersect. 
Let $X$ be the number of deleted edges.
There are fewer than $n^5$ pairs of edges that intersect, 
so $\mb{E}X < c^2 n/16$, so $\mb{P}(X < c^2 n/8) \ge 1/2$.
Also, the number of chosen absorbers $N_T$ for any triple $T$
is binomial with mean at least $c^2 n/4$, so whp all $N_T > c^2 n/8$.
Thus there is a choice of $A$ such that every $T$
has an absorber in $A$. This completes the proof
of the approximate version of \cite{RRS},
i.e.\ that minimum codegree $n/2 + cn$
guarantees a perfect matching. 

The idea for the exact result is to consider
an attempt to construct absorbers as above
under the weaker assumption $\dD(G) \ge n/2 - o(n)$.
It is not hard to see that absorbers exist unless
$G$ is close to one of the extremal examples.
The remainder of the proof (which we omit) is then
a stability analysis to show that the extremal
examples are locally optimal, and so optimal.

In the following two sections we will illustrate
two approaches to absorption for designs and hypergraph
decompositions, in the special case of triangle decompositions
of graphs, which is considerably simler, and so allows us to
briefly illustrate some (but not all) ideas 
needed for the general case. First we will conclude this
section by indicating why the basic method described above 
does not suffice.

Suppose we seek a triangle decomposition of a graph $G$
on $[n]$ with $e(G) = \OO(n^2)$ in which there  
is no space or divisibility obstruction:
we assume that $G$ is `tridivisible' 
(meaning that $3 \mid e(G)$ and all degrees are even)
and `triangle-regular' 
(meaning that there is a set $T$ of triangles in $G$
such that every edge is in $(1+o(1))tn$ triangles of $T$,
where $t>0$ and $n>n_0(t)$).
This is equivalent to a perfect matching in the auxiliary
$3$-graph $H$ with $V(H)=E(G)$ and 
$E(H) = \{\{ab,bc,ca\}: abc \in T\}$.
Note that $H$ is `sparse': we have $e(H) = O(v(H)^{3/2})$.
Triangle regularity implies that Pippenger's generalisation
(see \cite{PS}) of the R\"odl nibble can be applied
to give an almost perfect matching in $H$,
so the outstanding question is whether there is an absorber.

Let us consider a potential random construction
of an absorber $A$ in $H$. It will contain at most
$O(n^2)$ triangles, so the probability of any triangle
(assuming no heavy bias) will be $O(n^{-1})$.
On the other hand, to absorb some fixed (tridivisible) $S \sub E(G)$,
we need $A$ to contain a set $A_S$ of $a$ edge-disjoint triangles
(for some constant $a$) such that $S \cup A_S$ 
has a triangle decomposition $B_S$,
so we need $\OO(n^a)$ such $A_S$ in $G$.
To see that this is impossible, we imagine selecting the triangles
of $A_S$ one at a time and keeping track of the number $E_S$ of edges
that belong to a unique triangle of $S \cup A_S$.
If a triangle uses a vertex that has 
not been used previously then it increases $E_S$,
and otherwise it decreases $E_S$ by at most $3$.
We can assume that no triangle is used in both $A_S$ and $B_S$,
so we terminate with $E_S=0$. Thus there can be 
at most $3a/4$ steps in which $E_S$ increases,
so there are only $O(n^{3a/4})$ such $A_S$ in $G$.

The two ideas discussed below to overcoming this obstacle
can be briefly summarised as follows. In Randomised Algebraic
Construction (introduced in \cite{Kexist}), 
instead of choosing independent random triangles
for an absorber, they are randomly chosen according to a
superimposed algebraic structure that has `built-in' absorbers.
In Iterative Absorption (used for designs and
decompositions in \cite{GKLO,GKLO2}), 
instead of a single absorption step, 
there is a sequence of absorptions, each of which replaces 
the current subgraph of uncovered edges by an `easier' subgraph, 
until we obtain $S$ that is so simple that it can be absorbed 
by an `exclusive' absorber put aside at the beginning 
of the proof for the eventuality that we end up with $S$.
This is a powerful idea with many other applications
(see the survey \cite{KOicm}).

\section{Iterative Absorption}

Here we will sketch an application of iterative
absorption to finding a triangle decomposition 
of a graph $G$ with no space or divisibility
obstruction as in the previous subsection.
(We also make certain `extendability' assumptions 
that we will describe later when they are needed.)
Our sketch will be loosely based on a mixture 
of the methods used in \cite{BKLO} and \cite{GKLO},
thus illustrating some ideas of the general case
but omitting most of the technicalities.

The plan for the decomposition is to push the graph
down a `vortex', which consists of a nested sequence
$V(G) = V_0 \sups V_1 \sups \cdots \sups V_\tau$,
where $|V_i|=\tT|V_{i-1}|$ for each $i \in [\tau]$
with $n^{-1} \ll \tT \ll t$, and $|V_\tau|=O(1)$
(so $\tau$ is logarithmic in $n=v(G)$).
Suppose $G$ has a set $T$ of triangles such that
every edge is in $(1 \pm c)tn$ triangles of $T$,
where $n^{-1} \ll \tT \ll c, t$.
By choosing the $V_i$ randomly we can ensure that
each edge of $G[V_i]$ is in $(1 \pm 2c)t|V_i|$ 
triangles of $T_i = \{ f \in T: f \sub V_i \}$.
At step $i$ with $0 \le i \le t$ we will have
covered all edges of $G$ not contained in $V_i$
by edge-disjoint triangles, and also some edges
within $V_i$, in a suitably controlled manner,
so that we still have good triangle regularity in $G[V_i]$.

At the end of the process, the uncovered subgraph $L$
will be contained in $V_\tau$, so there are only
constantly many possibilities for $L$.
Before starting the process, for each tridivisible
subgraph $S$ of the complete graph on $V_\tau$
we put aside edge-disjoint `exclusive absorbers' $A_S$,
i.e.\ sets of edge-disjoint triangles in $G$
such that $S \cup A_S$ has a triangle decomposition $B_S$
(we omit here the details of this construction).
Then $L$ will be equal to one of these $S$, so replacing 
$A_S$ by $B_S$ completes the triangle decomposition of $G$.

Let us now consider the process of pushing $G$ down
the vortex; for simplicity of notation we describe 
the first step of covering all edges not within $V_1$.
The plan is to cover almost all of these edges by a nibble,
and then the remainder by a random greedy algorithm
(which will also use some edges within $V_1$).
At first sight this idea sounds suspicious,
as one would think that the triangle regularity 
parameter $c$ must increase substantially at each step,
and so the process could not be iterated
logarithmically many times before the parameters blow up.

However, quite suprisingly, if we make
the additional extendability assumption
that every edge is in at least $c' n^3$ copies of $K_5$
(where $c'$ is large compared with $c$ and $t \ge c'$),
then we can pass to a different
set of triangles which dramatically `boost' the regularity.
The idea (see \cite[Lemma 6.3]{GKLO}) is that
a relatively weak triangle regularity assumption
implies the existence of a perfect fractional triangle
decomposition, which can be interpreted as selection
probabilities (in the same spirit as \cite{KaLP})
for a new set of triangles that is much more regular.
A similar idea appears in the R\"odl-Schacht proof
of the hypergraph regularity lemma via
regular approximation (see \cite{RSc1}).
It may also be viewed as a `guided version'
of the self-correction that appears naturally
in random greedy algorithms (see \cite{BK2,FGM}).

Let us then consider $G^* = G \sm G[V_1] \sm H$,
where $H$ contains each edge of $G$ crossing
between $V_1$ and $V^* := V(G) \sm V_1$
independently with some small probability $p \ll c,\tT$.
(We reserve $H$ to help with the covering step.)
Then whp every edge of $G^*$ is in 
$(1 \pm c) tn \pm |V_1| \pm 3pn$
triangles of $T$ within $G^*$.
By boosting, we can find a set $T^*$ of triangles in $G^*$ 
such that every edge of $G^*$ is in $(1 \pm c_0) tn/2$
triangles of $T^*$, where $c_0 \ll p$.
By the nibble, we can choose a set of edge-disjoint
triangles in $T^*$ covering all of $G^*$
except for some `leave' $L$ of maximum degree $c_1 n$,
where we introduce new constants
$c_0 \ll c_1 \ll c_2 \ll p$.

Now we cover $L$ by two random greedy algorithms,
the first to cover all remaining edges in $V^*$
and the second to cover all remaining cross edges.
The analysis of these algorithms is not as difficult
as that of the nibble, as we have `plenty of space',
in that we only have to cover a sparse graph
within a much denser graph, whereas the nibble
seeks to cover almost all of a graph.
In particular, the behaviour of these algorithms
is well-approximated by a `binomial heuristic'
in which we imagine choosing random triangles
to cover the uncovered edges without worrying about
whether these triangles are edge-disjoint
(so we make independent choices for each edge).
In the actual algorithm we have to exclude any triangle
that uses an edge covered by a previous step of the 
algorithm, but if we are covering a sparse graph
one can show that whp at most half (say) of the choices
are forbidden at each step, so any whp estimate in
the binomial process will hold in the actual
process up to a factor of two.
(This idea gives a much simpler proof of the result of \cite{FHKS}.)

For the first greedy algorithm we consider 
each remaining edge in $V^*$ in some arbitrary order,
and when we consider $e$ we choose a triangle on $e$
whose two other edges are in $H$, and have not
been previously covered. In general 
we would make this choice uniformly at random, 
although the triangle case is sufficiently simple
that an arbitrary choice suffices; indeed,
there are whp at least $p^2 \tT n/2$ such triangles in $H$, 
of which at most $2c_1 n$ are forbidden due to using a 
previously covered edge (by the maximum degree of $L$). 
Thus the algorithm can be completed with arbitrary choices.

The second greedy algorithm for covering the
cross edges is more interesting 
(the analogous part of the proof for general designs
is the most difficult part of \cite{GKLO}). 
Let $H'$ denote the subgraph of cross edges
that are still uncovered. We consider each $x \in V^*$ 
sequentially and cover all edges of $H'$ incident to $x$
by the set of triangles obtained by adding $x$ to each
edge of a perfect matching $M_x$ in $G[H'(x)]$, i.e.\
the restriction of $G$ to the $H'$-neighbourhood of $x$.
We must choose $M_x$ edge-disjoint from $M_{x'}$ 
for all $x'$ preceding $x$, so an arbitrary choice
will not work; indeed, whp the degree of each vertex
$y$ in $G[H'(x)]$ is $(1 \pm c_2) p\tT tn$,
but our upper bound on the degree of $y$ in $H'$ 
may be no better than $pn$, so previous choices
of $M_{x'}$ could isolate $y$ in $G[H'(x)]$.

To circumvent this issue we choose random perfect matchings.
A uniformly random choice would work, but it is easier to
analyse the process where we fix many edge-disjoint matchings
in $G[H'(x)]$ and then choose one uniformly at random to be $M_x$.
We need some additional assumption to guarantee that
$G[H'(x)]$ has even one perfect matching 
(the approximate regularity only guarantees
an almost perfect matching).

One way to achieve this is to make
the additional mild extendability assumption
that every pair of vertices have at least $c' n$ 
common neighbours in $G[H'(x)]$, 
i.e.\ any adjacent pair of edges $xy, xy'$ in $G$
have at least $c' n$ choices of $z$
such that $xz$, $yz$ and $y'z$ are edges.
It is then not hard to see that a random
balanced bipartite subgraph of $G[H'(x)]$ whp
satisfies Hall's condition for a perfect matching.
Moreover, we can repeatedly delete $p^{3/2} \tT c'n$
perfect matchings in $G[H'(x)]$, as this maintains 
all degrees $(1 \pm 2\sqrt{p}) p\tT tn$ 
and codegrees at least $p\tT c'n/2$.

The punchline is that for any edge $e$ in $G[H'(x)]$
there are whp at most $2p^2 n$ earlier choices of $x'$
with $e$ in $G[H'(x')]$, and the random choice of $M_{x'}$
covers $e$ with probability at most $(p^{3/2} \tT c'n)^{-1}$,
so $e$ is covered with probability at most
$2p^2 n (p^{3/2} \tT c'n)^{-1} < p^{1/3}$, say.
Thus whp $G[H'(x)]$ still has sufficient
degree and codegree properties
to find the perfect matchings
described above, and the algorithm can be completed.
Moreover, any edge of $G[V_1]$ is covered with probability
at most $p^{1/3}$, so whp we maintain good triangle
regularity in $G[V_1]$ and can proceed down the vortex.

\section{Randomised Algebraic Construction}

Here we sketch an alternative proof (via our method of
Randomised Algebraic Construction from \cite{Kexist})
of the same result as in the previous subsection,
i.e.\ finding a triangle decomposition of a graph $G$ 
with certain extendability properties 
and no space or divisibility obstruction.
Our approach will be quite similar to that in \cite{Kcount},
except that we will illustrate the `cascade' approach to
absorption which is more useful for general designs.

As discussed above, we circumvent the difficulties
in the basic method for absorption by introducing
an algebraic structure with built-in absorbers.
Let $\pi:V(G) \to \mb{F}_{2^a} \sm \{0\}$ 
be a uniformly random injection, 
where $2^{a-2} < n \le 2^{a-1}$.
Our absorber (which in this context we call the `template')
is defined as the set $T$ of all triangles in $G$
such that $\pi(x)+\pi(y)+\pi(z)=0$.
Clearly $T$ consists of edge-disjoint triangles.
We let $G^* = \bigcup T$ be the underlying graph of
the template and suppress $\pi$, 
imagining $V(G)$ as a subset of $\mb{F}_{2^a}$.

Standard concentration arguments show that whp
$G \sm G^*$ has the necessary properties to apply
the nibble, so we can find a set $N$ of edge-disjoint
triangles in $G \sm G^*$ with leave
$L := (G \sm G^*) \sm \bigcup N$ 
of maximum degree $c_1 n$ 
(we use a similar hierarchy of very
small parameters $c_i$ as before).
To absorb $L$, it is convenient to first
`move the problem' into the template:
we apply a random greedy algorithm to cover $L$
by a set $M^c$ of edge-disjoint triangles,
each of which has one edge in $L$ and
the other two edges in $G^*$.
Thus some subgraph $S$ of $G^*$, 
which we call the `spill' has now 
been covered twice. The binomial heuristic
discussed in the previous subsection applies to
show that whp this algorithm is successful,
and moreover $S$ is suitably bounded.
(To be precise, we also ensure that each edge of $S$ 
belongs to a different triangle of $T$, and that
the union $S^*$ of all such triangles is $c_2$-bounded.)

The remaining task of the proof is to 
modify the current set of triangles
to eliminate the problem with the spill.
The overall plan is to find a `hole' 
in the template that exactly matches the spill.
This will consist of two sets of edge-disjoint triangles, 
namely $M^o$ (outer set) and $M^i$ (inner set),
such that $M^o \sub T$ and $\bigcup M^o$ is the 
disjoint union of $S$ and $\bigcup M^i$.
Then replacing $M^o$ by $M^i$ will fix the problem:
formally, our final triangle decomposition of $G$ 
is $M := N \cup M^c \cup (T \sm M^o) \cup M^i$.

We break down the task of finding the hole into several steps. 
The first is a refined form of the integral decomposition 
theorem of \cite{GJ,W4}, i.e.\ that there is an assignment
of integers to triangles so that total weight of triangles
on any edge $e$ is $1$ if $e \in S$ or $0$ otherwise. 
Our final hole can be viewed as such an assignment,
in which a triangle $f$ has weight $1$ if $f \in M^o$,
$-1$ if $f \in M^i$, or $0$ otherwise. We intend to start
from some assignment and repeatedly modify it by random greedy
algorithms until it has the properties required for the hole.
As discussed above, the success of such random greedy algorithms
requires control on the maximum degree, so our refined version
of \cite{GJ,W4} is that we can choose the weights $w_T$ 
on triangles with $\sum_{T: v \in T} |w_T| < c_3 n$
for every vertex $v$. (The proof is fairly simple, 
but the analogous statement for general hypergraphs 
seems to be much harder to prove.)
Note that in this step we allow the use of any triangle
in $K_n$ (the complete graph on $V(G)$),
without considering whether they belong to $G^*$:
`illegal' triangles will be eliminated later.

Let us now consider how to modify assignments
of weights to triangles so as to obtain a hole.
Our first step is to ignore the requirement $M^o \sub T$,
which makes our task much easier, 
as $T$ is a special set of only $O(n^2)$ triangles.
Thus we seek a signed decomposition of $S$ within $G^*$,
i.e.\ an assignment from $\{-1,0,1\}$ to each triangle 
of $G^*$ so that the total weight on 
any $e$ is $1$ if $e \in S$ or $0$ otherwise,
and every edge appears in at most 
one triangle of each sign.

To achieve this, we start from the simple observation that 
the graph of the octahedron has $8$ triangles, which can be split
into two groups of $4$, each forming a triangle decomposition.
For any copy of the octahedron in $K_n$ 
we can add $1$ to the triangles of one decomposition
and subtract $1$ from the triangles of the other
without affecting the total weight of triangles on any edge.
We can use this construction to repeatedly eliminate
`cancelling pairs', consisting of two triangles on
a common edge with opposite sign. (There is a preprocessing
step to ensure that each triangle to be eliminated 
can be assigned to a unique such pair.)
In particular, as edges not in $G^*$ have weight $0$,
this will eliminate all illegal triangles.
The boundedness condition facilitates a random greedy algorithm 
for choosing edge-disjoint octahedra for these eliminations,
which constructs the desired signed decomposition of $S$.

Now we remember that we wanted the outer triangle
decomposition $M^o$ to be contained in the template $T$.
Finally, the algebraic structure will come into play,
in absorbing the set $M^+$ of positive triangles
in the signed decomposition. To see how this can be achieved,
consider any positive triangle $xyz$,
recall that vertices are labelled 
by elements of $\mb{F}_{2^a} \sm \{0\}$,
and suppose first for simplicity that $xyz$
is `octahedral', meaning that $G^*$
contains the `associated octahedron' of $xyz$,
defined as the complete $3$-partite graph $O$ 
with parts $\{x,y+z\}$, $\{y,z+x\}$, $\{z,x+y\}$.
Then $xyz$ is a triangle of $O$, and we note that $O$
has a triangle decomposition consisting entirely 
of template triangles, namely $\{x,y,x+y\}$,
$\{y+z,y,z\}$, $\{x,z+x,z\}$ and $\{y+z,z+x,x+y\}$.
Thus we can `flip' $O$ (i.e.\ add and subtract the two 
triangle decompositions as before) to eliminate $xyz$ 
while only introducing positive triangles that are in $T$.

The approach taken in \cite{Kcount} was to ensure
in the signed decomposition that every positive triangle
is octahedral, with edge-disjoint associated octahedra,
so that all positive triangles can be absorbed as
indicated above without interfering with each other.
For general designs, it is more convenient to define 
a wider class of triangles (in general hypergraph cliques)
that can be absorbed by the following two step process,
which we call a `cascade'. Suppose that we want to
absorb some positive triangle $xyz$. We look for some
octahedron $O$ with parts $\{x,x'\}$, $\{y,y'\}$, $\{z,z'\}$
such that each of the $4$ triangles of the decomposition
not using $xyz$ is octahedral. We can flip the associated 
octahedra of these triangles so as to include them 
in the template, and now $O$ is decomposed by template
triangles, so can play the role of an associated octahedron
for $xyz$: we can flip it to absorb $xyz$.
The advantage of this approach is that whp
any non-template $xyz$ has many cascades,
so no extra property of the signed decomposition
is required to complete the proof.
In general, there are still some conditions
required for a clique have many cascades,
but these are not difficult to ensure 
in the signed decomposition.

\section{Concluding remarks}

There are many other questions of Design Theory
that can be reformulated as asking whether a certain
(sparse) hypergraph has a perfect matching.
This suggests the (vague) meta-question of formulating
and proving a general theorem on the existence of perfect
matchings in sparse `design-like' hypergraphs
(for some `natural' definition of `design-like'
that is sufficiently general to capture
a variety of problems in Design Theory).
One test for such a statement is that it should
capture all variant forms of the basic existence question,
such as general hypergraph decompositions (as in \cite{GKLO})
or resolvable designs (the general form of 
Kirkman's original `schoolgirl problem',
solved for graphs by Ray-Chaudhuri and Wilson \cite{RW}).

In \cite{K2} we generalised the existence of combinatorial designs
to the setting of subset sums in lattices with coordinates
indexed by labelled faces of simplicial complexes.
This general framework includes the problem of decomposing
hypergraphs with extra edge data, such as colours and orders, 
and so incorporates a wide range of variations on the basic design problem, 
notably Baranyai-type generalisations, such as resolvable hypergraph designs,
large sets of hypergraph designs and decompositions of designs by designs.
Our method also gives approximate counting results,
which is new for many structures whose existence was previously known,
such as high dimensional permutations or Sudoku squares.
For an exposition of these results and further applications, see \cite{KLovasz70}.

Could we be even more ambitious?
To focus the ideas, one well-known longstanding open
problem is Ryser's Conjecture \cite{ryser}
that every Latin square of odd order has a transversal.
(A generalised form of this conjecture by Stein \cite{stein}
was recently disproved by Pokrovskiy and Sudakov \cite{PoSu}.)
To see the connection with hypergraph matchings, 
we associate to any Latin square
a tripartite $3$-graph in which the parts correspond
to rows, columns and symbols, and each cell of the square 
corresponds to an edge consisting of its own
row, column and symbol. A perfect matching in this $3$-graph
is precisely a transversal of the Latin square.
However, there is no obvious common structure to the
various possible $3$-graphs that may arise in this way,
which presents a challenge 
to the absorbing methods described in this article,
and so to formulating a meta-theorem
that might apply to Ryser's Conjecture.
The best known lower bound of $n-O(\log^2n)$ 
on a partial transversal (by Hatami and Shor \cite{HS})
has a rather different proof. Another generalisation
of Ryser's Conjecture by Aharoni and Berger \cite{AB}
concerning rainbow matchings in properly coloured
multigraphs has recently motivated the development 
of various other methods for such problems
not discussed in this article (see e.g.\ \cite{grww,KY,Po}).

Recalling the theme of random matchings discussed
in the introduction, it is unsurprising that it is
hard to say much about random designs, but for certain
applications one can extract enough from the proof in \cite{Kexist},
e.g.\ to show that whp a random Steiner Triple System
has a perfect matching (Kwan \cite{Kw}) 
or that one can superimpose a constant number
of Steiner Systems to obtain a bounded 
codegree high-dimensional expander 
(Lubotzky, Luria and Rosenthal \cite{LLR}).
Does the nascent connection between hypergraph matchings 
and high-dimensional expanders go deeper?

We conclude by recalling two longstanding open problems
from the other end of the Design Theory spectrum,
concerning $q$-graphs with $q$ of order $\sqrt{n}$
(the maximum possible), as opposed to the setting $n>n_0(q)$
considered in this article (or even the methods of \cite{KLP}
which can allow $q$ to grow as a sufficiently small power of $n$).

\medskip

\nib{Hadamard's Conjecture.} (\cite{had}) \\
There is an $n \times n$ orthogonal matrix $H$ with all entries 
$\pm n^{-1/2}$ iff $n$ is $1$, $2$ or divisible by $4$?

\medskip

\nib{Projective Plane Prime Power Conjecture.} (folklore) \\
There is a Steiner system with parameters $(k^2+k+1,k+1,2)$
iff $k$ is a prime power?


\begin{thebibliography}{99}

\bibitem{A} R. Aharoni, 
Ryser's Conjecture for tripartite 3-graphs,
{\em Combinatorica} 21:1--4 (2001).

\bibitem{AB} R. Aharoni and E. Berger, 
Rainbow matchings in $r$-partite $r$-graphs, 
{\em Electron. J. Combin.} 16 (2009).

\bibitem{AH} R. Aharoni and P. Haxell,
Hall's theorem for hypergraphs,
{\em J. Graph Theory} 35:83--88 (2000).

\bibitem{ABHKP}
P. Allen, J. B\"ottcher, H. H\`an, Y. Kohayakawa and Y. Person,
Blow-up lemmas for sparse graphs, arxiv:1612.00622.

\bibitem{AKS} N. Alon, J. H. Kim and J. Spencer, 
Nearly perfect matchings in regular simple hypergraphs, 
{\em Israel J. Math.} 100:171--187 (1997).

\bibitem{AFS} A.~Asadpour, U.~Feige and A.~Saberi, Santa Claus meets hypergraph matchings, 
\emph{Proc.~APPROX-RANDOM} (2008), 10--20.

\bibitem{BKLO}
B. Barber, D. K\"uhn, A. Lo and D. Osthus,
Edge-decompositions of graphs with high minimum degree,
{\em Adv. Math.} 288:337--385 (2016).

\bibitem{BKLMO}
B. Barber, D. K\"uhn, A. Lo, R. Montgomery and D. Osthus,
Fractional clique decompositions of dense graphs and hypergraphs,
{\em  J. Combin. Theory Ser. B} 127:148--186 (2017).

\bibitem{BS} A. Barvinok and A. Samorodnitsky,
Computing the partition function for perfect matchings in a hypergraph,
{\em Combin. Probab. Comput.} 20:815--825 (2011).

\bibitem{BB} P. Bennett and T. Bohman,
A natural barrier in random greedy hypergraph matching,
arXiv:1210.3581.

\bibitem{Boh} T. Bohman, The triangle-free process, 
{\em Adv. Math.} 221:1653--1677 (2009).

\bibitem{BFL} T. Bohman, A. Frieze and E. Lubetzky,
Random triangle removal,
{\em Adv. Math.} 280:379--438 (2015).

\bibitem{BK} T. Bohman and P. Keevash, 
The early evolution of the H-free process, 
{\em Invent. Math.} 181:291--336 (2010). 

\bibitem{BK2} T. Bohman and P. Keevash, 
Dynamic concentration of the triangle-free process,
arXiv:1302.5963.

\bibitem{CD} C. J. Colbourn and J. H. Dinitz, 
{\em Handbook of Combinatorial Designs}, 
2nd ed. Chapman \& Hall / CRC, Boca Raton (2006).

\bibitem{CFMR} C. Cooper, A. Frieze, M. Molloy and B. Reed, 
Perfect matchings in random r-regular s-uniform hypergraphs, 
{\em Combin. Probab. Comput.} 5:1--14 (1996).

\bibitem{C} G. Cornuejols,
{\em Combinatorial Optimization: Packing and Covering},
Society for Industrial and Applied Mathematics (2001).

\bibitem{DK} P. Devlin and J. Kahn,
Perfect fractional matchings in k-out hypergraphs,
arXiv:1703.03513. 

\bibitem{DJPR}
E. Davies, M. Jenssen, W. Perkins and B. Roberts,
Tight bounds on the coefficients of partition functions via stability,
arXiv:1704.07784.

\bibitem{dross} F. Dross,
Fractional triangle decompositions 
in graphs with large minimum degree,
{\em SIAM J. Disc. Math.} 30:36--42 (2016).

\bibitem{dirac} G.~A.~Dirac, Some theorems on abstract graphs, 
\emph{Proc.~London Math.~Soc.} 2:69--81 (1952).

\bibitem{edmonds} J.~Edmonds, Paths, trees, and flowers, 
\emph{Canad.~J.~Math.} 17:449--467 (1965).

\bibitem{eg} J. Egerv\'ary, 
Matrixok kombinatorius tulajdons\'agair\'ol,
Matematikai \'es Fizikai Lapok 38:16--28 (1931).

\bibitem{e81} P. Erdo\H{o}s, 
On the combinatorial problems which 
I would most like to see solved, 
{\em Combinatorica} 1:25--42 (1981).

\bibitem{ematch} P. Erdo\H{o}s, 
A problem on independent r-tuples, 
{\em Ann. Univ. Sci. Budapest}, 8:93--95 (1965).

\bibitem{EH} P. Erd\H{o}s and H. Hanani, 
On a limit theorem in combinatorial analysis, 
{\em Publicationes Mathematicae Debrecen} 10:10--13 (1963).

\bibitem{FHKS} A. Ferber, R. Hod, M. Krivelevich and B. Sudakov, 
A construction of almost Steiner systems, 
{\em J. Combin. Designs} 22:488--494 (2014).

\bibitem{FGM} G. Fiz Pontiveros, S. Griffiths and R. Morris,
The triangle-free process and R(3,k), arXiv:1302.6279.

\bibitem{FR} P. Frankl and V. R\"odl, 
Near perfect coverings in graphs and hypergraphs, 
{\em Europ. J. Combin.} 6:317--326 (1985).

\bibitem{FT} P. Frankl, N. Tokushige,
Invitation to intersection problems for finite sets,
{\em J. Combin. Theory Ser. A} 144:157--211 (2016).

\bibitem{FKM} S. Friedland, E. Krop, and K. Markstr\"om,
On the number of matchings in regular graphs,
{\em Electronic J. Combin.} 15:R110 (2008).

\bibitem{fur} Z. F\"uredi, 
Matchings and covers in hypergraphs, 
{\em Graphs Combin.} 4:115--206 (1988).

\bibitem{GS} D. Gale and L.S. Shapley,
College admissions and the stability of marriage,
{\em Amer. Math. Monthly} 69:9--14 (1962).

\bibitem{grww} P. Gao, R. Ramadurai, I. Wanless and N. Wormald,
Full rainbow matchings in graphs and hypergraphs, arXiv:1709.02665.

\bibitem{GKLO} S. Glock, D. K\"uhn, A. Lo and D. Osthus,
The existence of designs via iterative absorption,
arXiv:1611.06827.

\bibitem{GKLO2} S. Glock, D. K\"uhn, A. Lo and D. Osthus,
Hypergraph $F$-designs for arbitrary $F$,
arXiv:1706.01800.

\bibitem{Gow} W. T. Gowers,
Hypergraph Regularity and the multidimensional Szemer\'edi Theorem, 
{\em Annals of Math.} 166:897--946 (2007).

\bibitem{Gr} D. A. Grable, More-than-nearly perfect packings and partial designs, 
{\em Combinatorica} 19:221--239 (1999).

\bibitem{GJ} J. E. Graver and W. B. Jurkat, 
The module structure of integral designs, 
{\em J. Combin. Theory Ser. A} 15:75--90 (1973).

\bibitem{had} J. Hadamard, 
R\'esolution d'une question relative aux d\'eterminants, 
{\em Bull. des Sciences Math.} 17:240--246 (1893). 

\bibitem{hall} P. Hall, On Representatives of Subsets, 
{\em J. London Math. Soc.} 10:26--30 (1935).

\bibitem{han} J. Han,
Decision problem for Perfect Matchings in Dense k-uniform Hypergraphs, 
{\em Trans. Amer. Math. Soc.} 369:5197--5218 (2017). 

\bibitem{han2} J. Han, 
Near perfect matchings in k-uniform hypergraphs,
\emph{Combin. Probab. Comput.} 24:723--732 (2015).

\bibitem{HS} P. Hatami and P. W. Shor,
A lower bound for the length of a partial transversal in a
Latin square, {\em J. Combin. Theory Ser. A} 115:1103--1113 (2008).

\bibitem{haxsurvey} P. Haxell,
Independent transversals and hypergraph matchings - an elementary approach,
in: {\em Recent Trends in Combinatorics}, Springer (2016).

\bibitem{HR} P. E. Haxell and V. R\"odl,
Integer and fractional packings in dense graphs. 
{\em Combinatorica} 21:13--38 (2001).

\bibitem{hen} J. R. Henderson, 
Permutation Decomposition 
of (0-1)-Matrices and Decomposition Transversals, 
Ph.D. thesis, Caltech (1971).

\bibitem{J} M. Jerrum,
{\em Counting, Sampling and Integrating: Algorithms and Complexity},
Springer Science \& Business Media (2003).

\bibitem{JKV} A. Johansson, J. Kahn and V. Vu,
Factors in Random Graphs,
{\em Random Struct. Alg.} 33:1--28 (2008).

\bibitem{Ka} J. Kahn, Asymptotically good list-colorings, 
{\em J. Combin. Theory Ser. A} 73:1--59 (1996).

\bibitem{KaLP} J. Kahn, 
A linear programming perspective on the
Frankl-R\"odl-Pippenger Theorem,
{\em Random Struct. Alg.} 8:149--157 (1996).

\bibitem{Ka2} J. Kahn, A normal law for matchings,
{\em Combinatorica} 20:339-391 (2000).

\bibitem{KaKa} J. Kahn and M. Kayll,
On the stochastic independence properties of hardcore distributions,
{\em Combinatorica} 17:369-391 (1997).

\bibitem{Karp} R.~M.~Karp, Reducibility among combinatorial problems, 
\emph{Complexity of Computer Computations} (1972), 85--103. 

\bibitem{Kexist} P. Keevash, The existence of designs, arXiv:1401.3665.

\bibitem{K2} P. Keevash, The existence of designs II, arXiv:1802.05900.

\bibitem{KLovasz70} P. Keevash, Coloured and directed designs, preprint.

\bibitem{Kcount} P. Keevash, Counting designs, 
to appear in {\em J. Eur. Math. Soc.}

\bibitem{Kblowup}
P. Keevash, A hypergraph blowup lemma,
{\em Random Struct. Alg.} 39:275--376 (2011). 

\bibitem{Kturan}  
P. Keevash, Hypergraph Tur\'an Problems, {\em Surveys in Combinatorics},
Cambridge University Press, 83--140 (2011). 

\bibitem{KKM} P.~Keevash, F.~Knox and R.~Mycroft, 
Polynomial-time perfect matchings in dense hypergraphs, 
{\em Adv. Math.} 269:265--334 (2015).
A preliminary version appeared in {\em Proc. 45th ACM STOC} (2013).


\bibitem{KM} P.~Keevash and R.~Mycroft, 
A geometric theory for hypergraph matchings, 
{\em Mem. Amer. Math. Soc.} 233 (2014), monograph 1098.

\bibitem{KM2} P. Keevash and R. Mycroft, 
A multipartite Hajnal-Szemeredi theorem, 
{\em J. Combin. Theory Ser. B} 114:187--236 (2015). 

\bibitem{KY} P. Keevash and L. Yepremyan,
Rainbow matchings in properly-coloured multigraphs,
arXiv:1710.03041.

\bibitem{kenyon} R. Kenyon, The dimer model,
in: {\em Exact Methods in Low-dimensional Statistical Physics and Quantum Computing},
Oxford University Press (2010).

\bibitem{Kim} J. H. Kim, Nearly optimal partial Steiner systems,
{\em Electronic Notes in Discrete Mathematics} 7:74--77.

\bibitem{KKOT}
J. Kim, D. K\"uhn, D. Osthus and M. Tyomkyn,
A blow-up lemma for approximate decompositions, 
to appear in {\em Trans. Amer. Math. Soc.}


\bibitem{KSS} J. Koml\'os, G. N. S\'ark\"ozy and E. Szemer\'edi,
Blow-up lemma, {\em Combinatorica} 17:109--123 (1997).

\bibitem{KR} A. Kostochka and V. R\"odl, 
Partial Steiner systems and matchings in hypergraphs, 
{\em Random Struct. Alg.} 13:335--347 (1997).

\bibitem{konig} D. K\"onig,
Gr\'afok \'es m\'atrixok,
{\em Matematikai \'es Fizikai Lapok}, 38:116--119 (1931).

\bibitem{kriv} M. Krivelevich, 
Perfect fractional matchings in random hypergraphs, 
{\em Random Struct. Alg.} 9:317--334 (1996).

\bibitem{kuhn} H. W. Kuhn, 
The Hungarian Method for the assignment problem, 
{\em Naval Research Logistics Quarterly}, 2:83--97 (1955).

\bibitem{KOicm} D. K\"uhn and D. Osthus,
Hamilton cycles in graphs and hypergraphs: an extremal perspective, 
{\em Proc. ICM 2014} 4:381--406, Seoul, Korea (2014).

\bibitem{KOlarge} D. K\"uhn and D. Osthus,
Embedding large subgraphs into dense graphs,
{\em Surveys in Combinatorics}, 
Cambridge University Press, 137--167 (2009).

\bibitem{KLP} G. Kuperberg, S. Lovett and R. Peled, 
Probabilistic existence of regular combinatorial objects, 
{\em Geom. Funct. Anal.} 27:919--972 (2017).
Preliminary version in {\em Proc. 44th ACM STOC} (2012).

\bibitem{Kuz} N. Kuzjurin, 
On the difference between asymptotically good parkings and coverings, 
{\em European J. Combin.} 16:35--40 (1995).

\bibitem{Kw} M. Kwan,
Almost all Steiner triple systems have perfect matchings,
arXiv:1611.02246.

\bibitem{LRS} L. C. Lau, R. Ravi and M. Singh, 
{\em Iterative Methods in Combinatorial Optimization},
Cambridge University Press (2011). 

\bibitem{L} L. Lov\'asz,
A kombinatorika minimax t\'eteleir\"ol 
(On minimax theorems in combinatorics),
{\em Matematikai Lapok} 26:209--264 (1975). 

\bibitem{LP} L. Lov\'asz and M. Plummer,
{\em Matching Theory}, American Mathematical Society (2009).

\bibitem{LLR}
A. Lubotzky, Z. Luria and R. Rosenthal,
Random Steiner systems and bounded degree 
coboundary expanders of every dimension,
arXiv:1512.08331.

\bibitem{luria} Z. Luria, 
New bounds on the number of n-queens configurations,
arXiv:1705.05225.

\bibitem{NRS} B. Nagle, V. R\"odl and M. Schacht,
The counting lemma for regular $k$-uniform hypergraphs,
{\em Random Struct. Alg.} 28:113--179 (2006).

\bibitem{NW} C. St. J. A. Nash-Williams, 
An unsolved problem concerning decomposition of graphs into triangles, 
{\em Combinatorial Theory and its Applications III}, North Holland (1970), 1179--1183.

\bibitem{PS} N. Pippenger and J. H. Spencer, 
Asymptotic behaviour of the chromatic index for hypergraphs, 
{\em J. Combin. Theory Ser. A} 51:24--42 (1989).

\bibitem{Po} A. Pokrovskiy, 
An approximate version of a conjecture of Aharoni and Berger, 
arXiv:1609.06346.

\bibitem{PoSu} A. Pokrovskiy and B. Sudakov, 
A counterexample to Stein's Equi-$n$-square Conjecture,
arXiv:1711.00429.

\bibitem{rad} Jaikumar Radhakrishnan, 
An entropy proof of Bregman's theorem, 
{\em J. Combin. Theory Ser. A} 77:80--83 (1997).

\bibitem{RW} D.K. Ray-Chaudhuri and R.M. Wilson,
Solution of Kirkman's schoolgirl problem, 
{\em Proc. Sympos. Pure Math.},
American Mathematical Society, XIX:187--203 (1971).

\bibitem{R} V. R\"odl, On a packing and covering problem,
{\em Europ. J. Combin.} 6:69--78 (1985).

\bibitem{RR} V.~R\"odl and A.~Ruci\'nski, Dirac-type questions for hypergraphs 
--- a survey (or more problems for Endre to solve),
\emph{An Irregular Mind (Szemer\'edi is 70)} 21:1--30 (2010).

\bibitem{RRS} V.~R\"odl, A.~Ruci\'nski and E.~Szemer\'edi, 
Perfect matchings in large uniform hypergraphs with large minimum collective degree, 
\emph{J.~Combin.~Theory Ser.~A} 113:613--636 (2009).

\bibitem{RSc1} V. R\"odl and M. Schacht, Regular partitions of hypergraphs: regularity lemma,
{\em Combin. Probab. Comput.} 16:833--885 (2007).

\bibitem{RSc2} V. R\"odl and M. Schacht,
Regular partitions of hypergraphs: counting lemmas,
{\em Combin. Probab. Comput.} 16:887--901 (2007).

\bibitem{RSST}
V. R\"odl and M. Schacht, M. H. Siggers and N. Tokushige,
Integer and fractional packings of hypergraphs,
{\em J. Combin. Theory Ser. B} 97:245--268 (2007).

\bibitem{RSk} V. R\"odl and J. Skokan, Regularity lemma for uniform hypergraphs,
{\em Random Struct. Alg.} 25:1--42 (2004).

\bibitem{ryser} H. Ryser, Neuere Probleme in der Kombinatorik,
Vortrage \"uber Kombinatorik, Oberwolfach, 69--91 (1967). 

\bibitem{schrijver} A. Schrijver,
{\em Combinatorial Optimization: Polyhedra and Efficiency},
Springer-Verlag (2003).

\bibitem{S} J. Spencer, Asymptotic packing via a branching process, 
{\em Random Struct. Alg.} 7:167--172 (1995).

\bibitem{stein} S. K. Stein, 
Transversals of Latin squares and their generalizations,
{\em Pacific J. Math.} 59:567-575 (1975).

\bibitem{St} E. Steinitz, 
\"Uber die Konstruction der Configurationen n (sub 3), 
Ph.D. thesis, Universit\"at Breslau (1894).

\bibitem{Sz} E. Szemer\'edi, Regular partitions of graphs,
Probl\`emes combinatoires et th\'eorie des graphes, 
{\em Colloq. Internat. CNRS} 260:399--401 (1978).

\bibitem{T} L. Teirlinck, Non-trivial t-designs without repeated blocks exist for all t, 
{\em Discrete Math.} 65:301--311 (1987).

\bibitem{Vu} V. Vu,
New bounds on nearly perfect matchings in hypergraphs: higher codegrees do help,
{\em Random Struct. Alg.} 17:29--63 (2000).

\bibitem{wachs} M. Wachs,
Topology of Matching, Chessboard and General Bounded Degree Graph Complexes,
{\em Algebra Univ.} 49:345--385 (2003).

\bibitem{WS} D. P. Williamson and D. B. Shmoys,
{\em The Design of Approximation Algorithms},
Cambridge University Press (2011).

\bibitem{RobinW} R. Wilson,
The early history of block designs,
{\em Rend. del Sem. Mat. di Messina} 9:267--276 (2003).

\bibitem{W1} R. M. Wilson, 
An existence theory for pairwise balanced designs 
I. Composition theorems and morphisms,
{\em J. Combin. Theory Ser. A} 13:220--245 (1972).

\bibitem{W2} R. M. Wilson, 
An existence theory for pairwise balanced designs 
II. The structure of PBD-closed sets and the existence conjectures,
{\em J. Combin. Theory Ser. A} 13:246--273 (1972).

\bibitem{W3} R. M. Wilson, 
An existence theory for pairwise balanced designs 
III. Proof of the existence conjectures,
{\em J. Combin. Theory Ser. A} 18:71--79 (1975).

\bibitem{W4} R. M. Wilson, 
The necessary conditions for t-designs are sufficient for something, 
{\em Utilitas Math.} 4:207--215 (1973).

\bibitem{W5} R. M. Wilson,
Signed hypergraph designs and diagonal forms for some incidence matrices,
{\em Des. Codes Cryptogr.} 17:289--297 (1999).

\bibitem{W6} R. M. Wilson,
Nonisomorphic Steiner Triple Systems,
{\em Math. Zeit.} 135:303--313 (1974).

\bibitem{Y} R. Yuster,
Combinatorial and computational aspects 
of graph packing and graph decomposition, 
{\em Computer Science Review} 1:12--26 (2007).

\bibitem{yi} Yi Zhao, 
Recent advances on Dirac-type problems for hypergraphs, 
in: {\em Recent Trends in Combinatorics}, Springer (2016).

\bibitem{yufei} Yufei Zhao, 
Extremal regular graphs: independent sets and graph homomorphisms,
arXiv:1610.09210. 



\end{thebibliography}
\end{document}